\documentclass[11pt, journal, draftclsnofoot, onecolumn]{IEEEtran}
\ifCLASSOPTIONcompsoc
\else
\fi

\ifCLASSINFOpdf
\else
\fi
\usepackage[parfill]{parskip}    
\usepackage{amsmath,amssymb,latexsym,enumerate, mathrsfs}
\usepackage{graphicx}
\usepackage{graphicx}
\usepackage{array}
\usepackage{times}
\usepackage{hyperref}
\DeclareMathOperator*{\argmin}{arg\,min}
\DeclareMathOperator*{\argmax}{arg\,max}

\newtheorem{lemma}{\bf Lemma}
\newtheorem{proposition}{\bf Proposition}

\newtheorem{remark}{Remark}

\newtheorem{definition}{\bf Definition}

\begin{document}
%
\title{On the Problem of Minimum Asymptotic Exit Rate for Stochastically Perturbed Multi-Channel Dynamical Systems}

\author{Getachew~K.~Befekadu,~\IEEEmembership{} and~Panos~J.~Antsaklis,~\IEEEmembership{}
\thanks{This work was supported in part by the National Science Foundation under Grant No. CNS-1035655. The first author acknowledges support from the Department of Electrical Engineering, University of Notre Dame.}
\IEEEcompsocitemizethanks{\IEEEcompsocthanksitem G. K. Befekadu is with the Department
of Electrical Engineering, University of Notre Dame, Notre Dame, IN 46556, USA.\protect\\
E-mail: gbefekadu1@nd.edu
\IEEEcompsocthanksitem P. J. Antsaklis is with the Department
of Electrical Engineering, University of Notre Dame, Notre Dame, IN 46556, USA.\protect\\
E-mail: antsaklis.1@nd.edu}}

\markboth{}%
{Shell \MakeLowercase{\textit{et al.}}: Bare Advanced Demo of IEEEtran.cls for Journals}
\IEEEcompsoctitleabstractindextext{%
\begin{abstract}
We consider the problem of minimizing the asymptotic exit rate with which the controlled-diffusion process of a stochastically perturbed multi-channel dynamical system exits from a given bounded open domain. In particular, for a class of admissible bounded linear feedback operators, we establish a connection between the asymptotic exit rate with which such a controlled-diffusion process exits from the given domain and the asymptotic behavior (i.e., a probabilistic characterization) of the principal eigenvalue of the infinitesimal generator, which corresponds to the stochastically perturbed dynamical system, with zero boundary conditions on the given domain. Finally, we briefly remark on the implication of our result for evaluating the performance of the associated deterministic multi-channel dynamical system, when such a dynamical system is composed with a set of (sub)-optimal admissible linear feedback operators.
\end{abstract}

}

\maketitle

\IEEEdisplaynotcompsoctitleabstractindextext

%
\IEEEpeerreviewmaketitle

\section{Introduction}	 \label{S1}
In this brief paper, we consider the problem of minimizing the asymptotic exit rate with which the controlled-diffusion process $x^{\epsilon}(t)$ exits from a given bounded open domain pertaining to the following stochastically perturbed multi-channel dynamical system
\begin{align}
d x^{\epsilon}(t) = A x^{\epsilon}(t) dt + \sum\nolimits_{i=1}^m B_i u_i(t) dt + \sqrt{\epsilon} \sigma(x^{\epsilon}(t))dW(t), \,\, x^{\epsilon}(0) = x_0, \label{Eq1}
\end{align}
where
\begin{itemize}
\item[-] $A \in \mathbb{R}^{d \times d}$, $B_i \in \mathbb{R}^{d \times r_i}$, $\epsilon$ is a small positive number (which represents the level of random perturbation in the system),
\item[-] $\sigma \colon \mathbb{R}^{d} \rightarrow \mathbb{R}^{d \times d}$ is Lipschitz with the least eigenvalue of $\sigma(\cdot)\sigma^T(\cdot)$ uniformly bounded away from zero, i.e., 
\begin{align*}
 \sigma(x)\sigma^T(x)  \ge \kappa I_{d \times d} , \quad \forall x \in \mathbb{R}^{d},
\end{align*}
for some $\kappa > 0$,
\item[-] $W(\cdot)$ is a $d$-dimensional standard Wiener process,
\item[-] $x^{\epsilon}(\cdot) \in \mathcal{X} \subseteq \mathbb{R}^{d}$ is the state trajectory of the system,
\item[-] {$u_i(\cdot)$ is a $\,\mathcal{U}_i$-valued measurable control process to the $i$th-channel (i.e., an admissible control from the measurable set $\mathcal{U}_i\subset \mathbb{R}^{r_i}$) such that for all $t > s$, $W(t)-W(s)$ is independent of $u_i(\nu)$ for $\nu \le s$ and }
\begin{align*}
\mathbb{E} \int_{0}^{t_1} \vert u(t)\vert^2 dt < \infty, \quad \forall t_1 \ge 0,
\end{align*}
{where $u(\cdot) \triangleq (u_1(\cdot), u_2(\cdot), \ldots, u_m(\cdot)) \in \prod\nolimits_{i=1}^m \mathcal{U}_i$.}
\end{itemize}
 
Let $D \subset \mathbb{R}^{d}$ be a bounded open domain with smooth boundary (i.e., $\partial D$ is a manifold of class $C^2$). Moreover, denote by $C_{0T}([0,T], \mathbb{R}^d)$ the space of all continuous functions $\varphi(t)$, $t \in [0,\, T]$, with range in $\mathbb{R}^d$; and, in this space, we define the following metric  
\begin{align}
\rho_{0T}(\varphi, \psi) =\sup_{t \in [0,\, T]} \Bigl\vert \varphi(t) - \psi(t) \Bigr\vert, \label{Eq2}
\end{align}
when $\varphi(t)$, $\psi(t)$ belong to $C_{0T}([0,T], \mathbb{R}^d)$.  If $\Phi$ is a subset of the space $C_{0T}([0,T], \mathbb{R}^d)$, then we define
\begin{align}
d_{0T}(\psi, \Phi) =\sup_{\varphi \in \Phi} \rho_{0T}(\psi(t),\varphi(t)). \label{Eq3}
\end{align}

{In what follows, we consider a particular class of admissible controls $u_i(\cdot) \in \mathcal{U}_i$ of the form $u_i(t)=\bigl(\mathcal{K}_i x^{\epsilon}\bigr)(t)$, $\forall t \ge 0$, where $\mathcal{K}_i$, for $i =1, 2, \ldots, m$, is a real, continuous $r_i \times d$ matrix function such that}
\begin{align}
  \mathscr{K} \subseteq \Biggl\{\underbrace{\bigl(\mathcal{K}_1, \mathcal{K}_2, \ldots, \mathcal{K}_m\bigr)}_{\substack{\triangleq \mathcal{K}}} \in \prod\nolimits_{i=1}^m \mathscr{K}_i[\mathcal{X}, \mathcal{U}_i] \biggm\lvert \phi\bigl(t; 0, x_0, (\mathcal{K}x^{0})(t)\bigr) \in \Omega, \notag \\ \forall t \ge 0, \, \, \forall x_0 \in \Omega \Biggr\}, \label{Eq4} 
\end{align}
{where $\mathscr{K}_i[\mathcal{X}, \mathcal{U}_i]$ is a closed subspace of bounded linear feedback operators from $\mathcal{X}$ to $\,\mathcal{U}_i$ and $\Omega$ is a bounded open set in $D \cup \partial D$ that contains the origin $0$. Moreover, $\phi\bigl(t; 0, x_0, (\mathcal{K} x^{0})(t)\bigr)$ is the unique solution for}
\begin{align}
  \dot{x}^{0}(t)=Ax^{0}(t) + \sum\nolimits_{i=1}^m B_i \bigl(\mathcal{K}_i x^{0}\bigr)(t),  \,\, x^{0}(0) = x_0 \in \Omega,  \label{Eq5} 
\end{align}
{that corresponds to the deterministic multi-channel dynamical system, when such a dynamical system is composed with a certain $m$-tuple of admissible linear feedback operators $\mathcal{K} \in \mathscr{K}$.}

Note that the infinitesimal generator pertaining to the controlled-diffusion process $x^{\epsilon}(t)$ of Equation~\eqref{Eq1}, when $u_i(t)=\bigl(\mathcal{K}_i x^{\epsilon}\bigr)(t)$, for $t \ge 0$ and $i =1, 2, \ldots, m$, is given by
\begin{align}
 \mathcal{L}_{\epsilon}^{\mathcal{K}}(\cdot)(x) = \Bigl \langle \bigtriangledown(\cdot), \Bigl(A x + \bigl(B, \mathcal{K}\bigr)x \Bigr) \Bigr\rangle + \frac{\epsilon}{2} \operatorname{tr}\Bigl \{\sigma(x)\sigma^T(x)\bigtriangledown^2(\cdot) \Bigr\}, \label{Eq6}
\end{align}
where $\bigl(B, \mathcal{K}\bigr)x(\cdot) = \sum\nolimits_{i=1}^m B_i \bigl(\mathcal{K}_i x\bigr)(\cdot)$ for all $t \ge 0$.

Let $\tau_D^{\epsilon}$ be the first exit-time for the controlled-diffusion process $x^{\epsilon}(t)$ from the domain $D$, i.e., 
\begin{align}
\tau_D^{\epsilon} = \inf \bigl\{ t > 0 \, \bigl\vert \, x^{\epsilon}(t) \in \partial D \bigr\}, \label{Eq7}
\end{align}
which also depends on the $m$-tuple of linear feedback operators $\mathcal{K} \in \mathscr{K}$ and, more precisely, on the behavior of the solutions to the deterministic dynamical system of Equation~\eqref{Eq5}. Moreover, let us denote by $\lambda_{\epsilon}^{\mathcal{K}}$ the principal eigenvalue of the infinitesimal generator $-\mathcal{L}_{\epsilon}^{\mathcal{K}}$ with zero boundary conditions on $\partial D$ which is given by
\begin{align}
\lambda_{\epsilon}^{\mathcal{K}} = - \limsup_{T \rightarrow \infty} \frac{1} {T} \log \mathbb{P}_{\epsilon}^{\mathcal{K}} \bigl\{\tau_D^{\epsilon} > T \bigr\}, \label{Eq8}
\end{align}
where the probability $\mathbb{P}_{\epsilon}^{\mathcal{K}}\bigl\{\cdot\bigr\}$ is conditioned on the initial point $x_0 \in D$ as well as on the class of admissible linear feedback operators $\mathscr{K}$.

{Next, let us introduce the following definition (i.e., the maximum closed invariant set for the deterministic dynamical system under the action of the class of linear feedback operators $\mathscr{K}$) which is useful in the following section.}
\begin{definition} \label{DFN1}
{A set $\Lambda_D^{\mathcal{K}} \subset D \cup \partial D$ is called the maximum closed invariant set for the deterministic dynamical system of Equation}~\eqref{Eq5} {(under the action of the $m$-tuple of linear feedback operators $\mathcal{K} \in \mathscr{K}$), if any set $\Omega \subset D \cup \partial D$, for some $\mathcal{K} \in \mathscr{K}$, satisfying the property}
\begin{align}
 \phi\bigl(t; 0, x_0, (\mathcal{K} x^{0})(t)\bigr) \in \Omega, \,\, \forall t \ge 0, \,\, \forall x_0 \in \Omega \label{Eq9}
\end{align}
{is a subset of $\Lambda_D^{\mathcal{K}}$.}
 \end{definition}

In Section~\ref{S2}, we provide an estimate for the asymptotic exit rate with which the controlled-diffusion process $x^{\epsilon}(t)$ exits from the domain $D$. In particular, we minimize the following quantity
\begin{align}
\lambda_{\epsilon}^{\mathcal{K}} = -\limsup_{T \rightarrow \infty} \frac{1}{T} \log \mathbb{P}_{\epsilon}^{\mathcal{K}} \bigl\{ \tau_{D}^{\epsilon} > T\bigr\}, \label{Eq10}
\end{align}
with respect to the admissible controls $u_i(\cdot) \in \mathcal{U}_i$ of the form $u_i(t)=\bigl(\mathcal{K}_i x^{\epsilon}\bigr)(t)$, for $t \ge 0$ and for all $i \in \{1,2, \ldots m\}$. Note that if the domain $D$ contains an equilibrium point for the deterministic dynamical system of Equation~\eqref{Eq5}, when such a dynamical system is composed with the $m$-tuple of linear feedback operators $\mathcal{K} \in \mathscr{K}$. Then, the principal eigenvalue $\lambda_{\epsilon}^{\mathcal{K}}$ of the infinitesimal generator $-\mathcal{L}_{\epsilon}^{\mathcal{K}}$ with zero boundary conditions on $\partial D$ is equal to $\lambda_{\epsilon}^{\mathcal{K}} = \epsilon^{-1} r(\mathcal{K}) + o(\epsilon^{-1})$ as $\epsilon \rightarrow 0$ (e.g., see \cite{VenFre70}, \cite{Ven72} or \cite{Day83}). {On the other hand, if the maximum closed invariant set for the deterministic dynamical system under the action of the $m$-tuple of linear feedback operators $\mathcal{K} \in \mathscr{K}$ is nonempty. Then, the following asymptotic condition also holds true}
\begin{align}
-\lim_{\epsilon \rightarrow 0} \, \limsup_{T \rightarrow \infty} \frac{1} {T} \log \mathbb{P}_{\epsilon}^{\mathcal{K}} \bigl\{\tau_D^{\epsilon} > T \bigr\} < \infty, \,\,x_0 \in D. \label{Eq11}
\end{align}

\begin{remark} \label{R1}
{Note that such an asymptotic behavior of $\frac{1}{T} \log \mathbb{P}_{\epsilon}^{\mathcal{K}} \bigl\{\tau_D^{\epsilon} > T \bigr\}$ as $\epsilon \rightarrow 0$ and $T  \rightarrow \infty$, determines whether the dynamical system of Equation}~\eqref{Eq5} {has a maximum closed invariant set in $D \cup \partial D$ or not} (see \cite[Theorem~2.1]{Kif81}).
\end{remark}

Moreover, the principal eigenvalue $\lambda_{\epsilon}^{\mathcal{K}}$ turns out to be the boundary value between those $R < r(\mathcal{K})$ for which $\mathbb{E}_{\epsilon}^{\mathcal{K}}\bigl\{\exp(\epsilon^{-1} R \tau_D^{\epsilon})\bigr\} < \infty$ and those $R > r(\mathcal{K})$ for which $\mathbb{E}_{\epsilon}^{\mathcal{K}}\bigl\{\exp(\epsilon^{-1} R\tau_D^{\epsilon})\bigr\} = \infty$, where $r(\mathcal{K})$ is given by the following
\begin{align}
r(\mathcal{K}) = \limsup_{T \rightarrow \infty} \,\inf_{\substack{\varphi(t) \in C_{0T}([0,T], \mathbb{R}^d)\\ \varphi(0)=x_0}} \frac{1}{T} \biggl \{ S_{0T}^{\mathcal{K}}(\varphi(t)) \, \Bigl \vert \, \varphi(t) \in D \cup \partial D,\,\, \forall t \in [0, T] \biggr\}. \label{Eq12}
\end{align}
In general, such an asymptotic analysis involves minimizing the following action functional
\begin{align}
 S_{0T}^{\mathcal{K}}(\varphi(t)) =  \frac{1}{2} \int_{0}^{T}\biggl\Vert\frac{d\varphi(t)}{dt} - \Bigl(A \varphi(t) + \bigl(B, \mathcal{K}\bigr)\varphi(t) \Bigr) \biggr\Vert^2dt, \label{Eq13}
\end{align}
where
\begin{align}
\biggl\Vert\frac{d\varphi(t)}{dt} -  \Bigl(A \varphi(t) + \bigl(B, \mathcal{K}\bigr)\varphi(t) \Bigr) \biggr\Vert^2 = \biggl[\frac{d\varphi(t)}{dt} - \Bigl(A \varphi(t) + \bigl(B, \mathcal{K}\bigr)\varphi(t) \Bigr) \biggr]^T \notag \\
 \times \Bigl(\sigma(\varphi(t))\sigma^T(\varphi(t))\Bigr)^{-1} \biggl[\frac{d\varphi(t)}{dt} -  \Bigl(A \varphi(t) + \bigl(B, \mathcal{K}\bigr)\varphi(t) \Bigr) \biggr], \label{Eq14}
\end{align}
with $\bigl(B, \mathcal{K}\bigr)\varphi(t) = \sum\nolimits_{i=1}^m B_i \bigl(\mathcal{K}_i \varphi\bigr)(t)$ and $\varphi(t) \in C_{0T}([0,T], \mathbb{R}^d)$ is absolutely continuous.

Note that estimating the asymptotic exit rate with which the controlled-diffusion process $x^{\epsilon}(t)$ exits from the domain $D$ is related to a singularly perturbed eigenvalue problem. For example, the asymptotic behavior for the principal eigenvalue corresponding to the following eigenvalue problem
\begin{align}
\left.\begin{array}{c}
 -\mathcal{L}_{\epsilon}^{\mathcal{K}}\,\upsilon_{(\epsilon, x_0)}^{\mathcal{K}} = \lambda_{\epsilon}^{\mathcal{K}}\,\upsilon_{(\epsilon, x_0)}^{\mathcal{K}} \quad \text{in} \quad D\\
\quad\quad ~ \upsilon_{(\epsilon, x_0)}^{\mathcal{K}} = 0 \quad \text{on} \quad \partial D\\
\end{array}\right\}, \label{Eq15}
\end{align}
where $\upsilon_{(\epsilon, x_0)}^{\mathcal{K}} \in W_{loc}^{2, p} \cap C(D \cup \partial D)$, for $p > 2$, with $\upsilon_{(\epsilon, x_0)}^{\mathcal{K}} > 0$ on $D$, has been well studied in the past (e.g., see \cite{Day83} or \cite{DevFr78} in the context of an asymptotic behavior for the principal eigenfunction; and see \cite{Day87} or \cite{BiBor09} in the context of an asymptotic behavior for the equilibrium density). {Specifically, the author in} \cite{Day87} {has also provided additional results in connection with the asymptotic behavior of the equilibrium density, when the latter (i.e., the asymptotic behavior of the equilibrium density) is associated with the boundary exit problem from the domain of attraction with an exponentially stable critical point for the stochastically perturbed dynamical system} (e.g., see \cite{She91} or \cite{DayDa84}). 

{Before concluding this section, it is worth mentioning that, some interesting studies on the asymptotic behavior of dynamical systems with small random perturbations have been reported in control theory literature (to mention a few, e.g., see} \cite{Fle85}, \cite{CotFoMa83} or \cite{Zab85} {in the context of stochastic control approach to large deviation problems -- based on the Ventcel-Freidlin estimates} \cite{VenFre70} (cf. \cite[Chapter~14]{Fre76} or \cite{FreWe84}); {and see also} \cite{Sas83} {in the context of jump phenomena in nonlinear dynamical systems).}

\section{Main Results} \label{S2}
In this section, we present our main result -- where we establish a connection between the asymptotic exit rate with which the controlled-diffusion process $x^{\epsilon}(t)$ exits from the domain $D$ and the asymptotic behavior of the principal eigenvalue of the infinitesimal generator $-\mathcal{L}_{\epsilon}^{\mathcal{K}}$ with zero boundary conditions on $\partial D$.

In what follows, we state the following lemmas that will be useful for proving our main results (see \cite[Theorem~1.1, Theorem~1.2 and Lemma~9.1]{VenFre70} or \cite{Ven73}; and see \cite[pp.\,332--340]{Fre76} for additional discussions).

\begin{lemma} \label{L1}
For any $\alpha >0$, $\delta > 0$ and $\gamma > 0$, there exists an $\epsilon_0 > 0$ such that
\begin{align}
\mathbb{P}_{\epsilon}^{\mathcal{K}}\Bigl\{\rho_{0T}\bigl(x^{\epsilon}(t), \varphi(t)\bigr) < \delta \Bigr\} \ge \exp \Bigl\{ - \epsilon^{-1}\bigl(S_{0T}^{\mathcal{K}}(\varphi(t)) + \gamma \bigr)\Bigl\}, \,\, \forall  \epsilon \in (0, \epsilon_0), \label{Eq16}
\end{align}
where $\varphi(t)$ is any function in $C_{0T}([0,T], \mathbb{R}^d)$ for which $S_{0T}^{\mathcal{K}}(\varphi(t)) < \alpha$ and $\varphi(0)=x_0$.
\end{lemma}

\begin{lemma} \label{L2}
For any $\alpha >0$, $\delta > 0$ and $\gamma > 0$, there exists an $\epsilon_0 > 0$ such that
\begin{align}
\mathbb{P}_{\epsilon}^{\mathcal{K}}\Bigl\{d_{0T}\bigl(x^{\epsilon}(t), \Phi_{x_0, \alpha} \bigr) \ge \delta \Bigr\} \le \exp \Bigl\{ - \epsilon^{-1}\bigl(\alpha - \gamma \bigr)\Bigl\}, \,\, \forall  \epsilon \in (0, \epsilon_0), \label{Eq17}
\end{align}
where
\begin{align}
\Phi_{x_0, \alpha} = \Bigl\{\varphi(t) \in C_{0T}([0,T], \mathbb{R}^d) \, \Bigl \vert \,\varphi(0) = x_0\,\, \text{and}  \,\, S_{0T}^{\mathcal{K}}(\varphi(t)) < \alpha \Bigl\}. \label{Eq18}
\end{align}
\end{lemma}

\begin{lemma} \label{L3}
Let $D_{+\delta}$ denote a $\delta$-neighborhood of $D$ and let $D_{-\delta}$ denote the set of points in $D$ at a distance greater than $\delta$ from the boundary $\partial D$. Then, for sufficiently small $\delta > 0$ and for any $\mathcal{K} \in \mathscr{K}$, the following estimates 
\begin{align}
\inf_{\substack{\varphi(t) \in C_{0T}([0,T], \mathbb{R}^d) \\ \varphi(0)=x_0}} \biggl \{ S_{0T}^{\mathcal{K}}(\varphi(t)) \, \Bigl \vert \,\varphi(t) \in D_{+\delta} \cup \partial D_{+\delta}, \,\, \forall t \in [0, T] \biggr \}, \label{Eq19}
\end{align}
and 
\begin{align}
\inf_{\substack{\varphi(t) \in C_{0T}([0,T], \mathbb{R}^d)\\ \varphi(0)=x_0}} \biggl \{ S_{0T}^{\mathcal{K}}(\varphi(t))\, \Bigl \vert \,\varphi(t) \in D_{-\delta} \cup \partial D_{-\delta}, \,\, \forall t \in [0, T] \biggr\}, \label{Eq20}
\end{align}
can be made arbitrarily close to each other. Furthermore, the same holds for
\begin{align}
\inf_{\substack{\varphi(t) \in C_{0T}([0,T], \mathbb{R}^d)\\ \varphi(0)=x,\,\varphi(T)=y}} \biggl \{ S_{0T}^{\mathcal{K}}(\varphi(t))\, \Bigl \vert \,\varphi(t) \in D_{\pm\delta} \cup \partial D_{\pm\delta}, \, \, \forall t \in [0, T] \biggr\}, \label{Eq21}
\end{align}
uniformly for any $x, y \in D_{-\delta}$.
\end{lemma}

The following proposition provides an estimate for the principal eigenvalue $\lambda_{\epsilon}^{\mathcal{K}}$ of the infinitesimal generator $-\mathcal{L}_{\epsilon}^{\mathcal{K}}$ with zero boundary conditions on $\partial D$. {Note that such an estimate for the principal eigenvalue is apparently related to the asymptotic exit rate with which the controlled-diffusion process $x^{\epsilon}(t)$ exits from the domain $D$ -- when the dynamical system of Equation~}\eqref{Eq1} {is composed with the $m$-tuple of linear feedback operators $\mathcal{K} \in \mathscr{K}$.}
\begin{proposition} \label{P1}
{Suppose that the maximum closed invariant set for the deterministic dynamical system of Equation}~\eqref{Eq5}, {under the action of the $m$-tuple of linear feedback operators $\mathcal{K} \in \mathscr{K}$, is nonempty.} Then, the principal eigenvalue $\lambda_{\epsilon}^{\mathcal{K}}$ of the infinitesimal generator $-\mathcal{L}_{\epsilon}^{\mathcal{K}}$ with zero boundary conditions on $\partial D$ satisfies 
\begin{align}
\lambda_{\epsilon}^{\mathcal{K}} = \epsilon^{-1} r(\mathcal{K}) + o(\epsilon^{-1}) \quad \text{as} \quad \epsilon \rightarrow 0, \label{Eq22}
\end{align}
where 
\begin{align}
r({\mathcal{K}}) = \limsup_{T \rightarrow \infty} \,\inf_{\substack{\varphi(t) \in C_{0T}([0,T], \mathbb{R}^d)\\ \varphi(0)=x_0}} \frac{1}{T} \biggl \{ S_{0T}^{\mathcal{K}}(\varphi(t)) \, \Bigl \vert \,\varphi(t) \in D \cup \partial D,\,\, \forall t \in [0, T] \biggr\}. \label{Eq23}
\end{align}
\end{proposition}

\begin{IEEEproof}
Suppose that $r({\mathcal{K}})$, for a certain $m$-tuple of linear feedback operators $\mathcal{K} \in \mathscr{K}$, exists.\footnote{Note that the existence of such a limit for $r({\mathcal{K}})$ can be easily established (e.g., see \cite{VenFre70}).} Then, using Lemma~\ref{L3}, one can show that $r({\mathcal{K}})$ also satisfies the following
\begin{align}
r({\mathcal{K}}) =  \sup_{x, y \in D} \Biggl\{ \limsup_{T \rightarrow \infty} \, \inf_{\substack{\varphi(t) \in C_{0T}([0,T], \mathbb{R}^d)\\ \varphi(0)=x,\,\varphi(T)=y}} \frac{1}{T} \biggl \{ S_{0T}^{\mathcal{K}}(\varphi(t))\, \Bigl \vert \,\varphi(t) \in D \cup \partial D,\,\, \forall t \in [0, T] \biggr\} \Biggr\}. \label{Eq24}
\end{align}
Next, let us show that, for sufficiently small $\epsilon > 0$, $\mathbb{E}_{x_0}^{\epsilon}\bigl\{\exp(\epsilon^{-1} R \tau_D^{\epsilon})\bigr\}$ tends to infinity, when $R > r({\mathcal{K}})$. If we choose a positive $\varkappa$ which is smaller than $(R-r({\mathcal{K}}))/3$ so that
\begin{align}
 \sup_{x, y \in D} \, \inf_{\substack{\varphi(t) \in C_{0T}([0,T], \mathbb{R}^d)\\ \varphi(0)=x,\,\varphi(T)=y}} \frac{1}{T} \biggl \{ S_{0T}^{\mathcal{K}}(\varphi(t))\, \Bigl \vert \,\varphi(t) \in D \cup \partial D,\,\, \forall t \in [0, T] \biggr\} < r({\mathcal{K}}) + \varkappa, \label{Eq25}
\end{align}
and, for sufficiently small $\delta > 0$,
\begin{align}
 \inf_{\substack{\varphi(t) \in C_{0T}([0,T], \mathbb{R}^d)\\ \varphi(0)=x,\,\varphi(T)=y}} \biggl \{ S_{0T}^{\mathcal{K}}(\varphi(t))\, \Bigl \vert \,\varphi(t) \in D_{-\delta} \cup \partial D_{-\delta}, \,\, \forall t \in [0, T] \biggr\} < T\bigl(r({\mathcal{K}}) + 2\varkappa\bigr), \label{Eq26}
\end{align}
for all $x, y \in D_{-\delta}$. Then, if we further let $\alpha = T\bigl(r({\mathcal{K}}) + 2 \varkappa\bigr)$ and $\gamma = T\varkappa$, from Lemma~\ref{L1}, there exits an $\epsilon_0 > 0$ such that
\begin{align}
 S_{0T}^{\mathcal{K}}(\varphi(t)) \le T\bigl(r({\mathcal{K}}) + 2 \varkappa\bigr), \quad \varphi(t) \in D_{-\delta} \cup \partial D_{-\delta}, \forall t \in [0, T], \label{Eq27}
\end{align}
for any $x, y \in D_{-\delta}$; and, moreover, we have the following probability estimate
\begin{align}
\mathbb{P}_{x}^{\epsilon} \Bigl\{ \tau_{D_{-\delta}}^{\epsilon} > T,\, x^{\epsilon}(\tau_{D_{-\delta}}) \in D_{-\delta} \Bigr\} & \ge \mathbb{P}_{x}^{\epsilon}\Bigl\{\rho_{0T}\bigl(x^{\epsilon}(t), \varphi(t)\bigr) < \delta \Bigr\}, \notag \\ 
&\ge \exp \bigl(-\epsilon^{-1}\bigl(S_{0T}^{\mathcal{K}}(\varphi(t)) + \gamma \bigr)\bigl), \notag \\
&\ge \exp \bigl( -\epsilon^{-1} T \bigl(r({\mathcal{K}}) + 3\varkappa \bigr)\bigl), \quad \forall  \epsilon \in (0, \epsilon_0), \label{Eq28}
\end{align}
where $\varphi(T) \in D_{-2\delta}$.

Let us define the following random events 
\begin{align}
\mathcal{A}_n = \Bigl\{ \tau_{D_{-\delta}}^{\epsilon} > nT,\, x^{\epsilon}(nT) \in D_{-\delta}\Bigl\}, \label{Eq29}
\end{align}
for $n \in \mathbb{N}_{+} \cup \{0\}$. Then, from the strong Markov property, we have
\begin{align}
\mathbb{P}_{x}^{\epsilon} \bigl\{ \mathcal{A}_n \bigl\} & \ge \mathbb{E}_{x}^{\epsilon} \chi_{\mathcal{A}_{n-1}} \mathbb{P}_{x_{(n-1)T}}^{\epsilon} \bigl\{ \mathcal{A}_1 \bigr\},\notag  \\
& \ge \mathbb{P}_{x}^{\epsilon} \bigl\{ \mathcal{A}_{n-1} \bigr\}  \inf_{y \in D_{-\delta}} \mathbb{P}_{y}^{\epsilon} \bigl\{ \mathcal{A}_1 \bigr\},\notag  \\
& \ge \exp \bigl( - \epsilon^{-1} n T\bigl(r({\mathcal{K}}) + 3\varkappa \bigr)\bigr) \quad \forall  \epsilon \in (0, \epsilon_0). \label{Eq30}
\end{align}
Note that, for an arbitrary $n$, we have the following
\begin{align}
\mathbb{E}_{x}^{\epsilon} \bigl\{ \exp \bigl( \epsilon^{-1} R \tau_{D_{-\delta}}^{\epsilon} \bigl) \bigl\} & \ge \exp \bigl( \epsilon^{-1} R nT \bigl) \mathbb{P}_{x}^{\epsilon}\bigl\{ \tau_{D_{-\delta}}^{\epsilon} >  n T \bigr\},\notag  \\
& \ge \exp \bigl( - \epsilon^{-1} n T\bigl(R - r({\mathcal{K}}) - 3\varkappa \bigr)\bigr), \quad \forall  \epsilon \in (0, \epsilon_0), \label{Eq31}
\end{align}
which tends to infinity as $n \rightarrow \infty$, i.e., $\mathbb{E}_{x}^{\epsilon} \bigl\{ \exp \bigl(\epsilon^{-1} R \tau_{D_{-\delta}}^{\epsilon} \bigl) \bigl\}=\infty$.

On the other hand, let us show that if $R < r({\mathcal{K}})$, then, for sufficiently small $\epsilon > 0$, $\mathbb{E}_{x}^{\epsilon} \bigl\{ \exp \bigl(\epsilon^{-1} R \tau_{D_{-\delta}}^{\epsilon} \bigl) \bigl\} \\ < \infty$. For $\varkappa < (R - r({\mathcal{K}}))/3$, let us choose $\delta$ so that
\begin{align}
\inf_{\substack{\varphi(t) \in C_{0T}([0,T], \mathbb{R}^d) \\ \varphi(0)=x_0}} \biggl \{ S_{0T}^{\mathcal{K}}(\varphi(t)) \, \Bigl \vert \,\varphi(t) \in D_{+\delta} \cup \partial D_{+\delta}, \, \, \forall t \in [0, T] \biggr \} > T\bigl(r({\mathcal{K}}) - 2 \varkappa \bigr). \label{Eq32}
\end{align}
From Lemma~\ref{L2}, with $\alpha = T\bigl(r({\mathcal{K}}) - 2 \varkappa \bigr)$ and $\gamma T\varkappa$, there exists an $\epsilon_0 > 0$ such that the distance between the set of functions $\psi(t)$, for $0 \le t \le T$, entirely lying in $D$ and any of the sets $\Phi_{x_0, \alpha}$ is at least a distance $\delta$; and, hence, we have the following probability estimate
\begin{align}
\mathbb{P}_{x}^{\epsilon} \Bigl\{ \tau_{D}^{\epsilon} > T \Bigr\} & \le \mathbb{P}_{x}^{\epsilon}\Bigl\{d_{0T}\bigl(x^{\epsilon}(t), \Phi_{x_0, \alpha}\bigr) \ge \delta \Bigr\},\notag \\ 
&\le \exp \bigl(-\epsilon^{-1} T \bigl(r({\mathcal{K}}) - 3\varkappa \bigr)\bigl), \quad \forall  \epsilon \in (0, \epsilon_0), \label{Eq33}
\end{align}
for any $x \in D$.

Then, using the Markov property, we have the following
\begin{align}
\mathbb{P}_{x}^{\epsilon} \Bigl\{ \tau_{D}^{\epsilon} > nT \Bigr\} &\le \exp \bigl( -\epsilon^{-1} nT \bigl(r({\mathcal{K}}) - 3\varkappa \bigr)\bigl), \quad \forall  \epsilon \in (0, \epsilon_0), \label{Eq34}
\end{align}
and
\begin{align}
\mathbb{E}_{x}^{\epsilon} \bigl\{ \exp \bigl( \epsilon^{-1} R \tau_{D}^{\epsilon} \bigl) \bigl\} & \le \sum\nolimits_{n=0}^{\infty} \exp \bigl( \epsilon^{-1} R (n +1) T \bigl) \mathbb{P}_{x}^{\epsilon}\bigl\{nT <  \tau_{D}^{\epsilon} \le (n+1)T \bigr\},\notag  \\
& \le \sum\nolimits_{n=0}^{\infty} \exp \bigl( \epsilon^{-1} R (n +1) T \bigl) \mathbb{P}_{x}^{\epsilon}\bigl\{ \tau_{D}^{\epsilon} > nT \bigr\},\notag  \\
& \le \sum\nolimits_{n=0}^{\infty} \exp \bigl( \epsilon^{-1} R T \bigl) \exp \bigl( -\epsilon^{-1} n T\bigl(R - r({\mathcal{K}}) + 3\varkappa \bigr)\bigr), \,\, \forall \epsilon \in (0, \epsilon_0), \label{Eq35}
\end{align}
which converges to a finite value, i.e., $\mathbb{E}_{x}^{\epsilon} \bigl\{ \exp \bigl( \epsilon^{-1} R \tau_{D}^{\epsilon} \bigl) \bigl\} < \infty$. Hence, $r({\mathcal{K}})$ is a boundary for which $\mathbb{E}_{x}^{\epsilon}\bigl\{\exp(\epsilon^{-1} r({\mathcal{K}}) \tau_D^{\epsilon})\bigr\}$ is finite. Then, from Equation~\eqref{Eq33} (cf. Equation\eqref{Eq28}), we have
\begin{align}
-\frac{1}{T} \log \mathbb{P}_{x}^{\epsilon} \Bigl\{ \tau_{D}^{\epsilon} > T \Bigr\} \le \epsilon^{-1}\bigl(r({\mathcal{K}}) - 3\varkappa \bigr), \quad \forall  \epsilon \in (0, \epsilon_0), \label{Eq36}
\end{align}
for any $x \in D$, where the left side tends to the principal eigenvalue $\lambda_{\epsilon}^{\mathcal{K}}$ as $T \rightarrow \infty$. This completes the proof.
\end{IEEEproof}

\begin{remark} \label{R2}
{The above proposition states that the principal eigenvalue $\lambda_{\epsilon}^{\mathcal{K}}$ of the infinitesimal generator $-\mathcal{L}_{\epsilon}^{\mathcal{K}}$ with zero boundary conditions on $\partial D$ is equal to $\lambda_{\epsilon}^{\mathcal{K}} = \epsilon^{-1} r(\mathcal{K}) + o(\epsilon^{-1})$ as $\epsilon \rightarrow 0$} (cf. Equation~\eqref{Eq36} and Equation~\eqref{Eq34}). {Furthermore, the principal eigenvalue is exactly equal to the boundary value between those positive $R < r({\mathcal{K}})$ for which $\mathbb{E}_{x}^{\epsilon}\bigl\{\exp(\epsilon^{-1} R \tau_D^{\epsilon})\bigr\} < \infty$ and those $R > r({\mathcal{K}})$  for which $\mathbb{E}_{x}^{\epsilon}\bigl\{\exp(\epsilon^{-1} R\tau_D^{\epsilon})\bigr\} = \infty$ for any $x \in D$.}
\end{remark}

\section{Further Remarks} \label{S3}
{In this section, we further comment on the implication of our main result -- when one is also interested in either evaluating the performance of a certain admissible linear feedback operators $\hat{\mathcal{K}} \in \mathscr{K}$ (or finding a set of (sub)-optimal admissible linear feedback operators $\bigl\{\mathcal{K}^{\ast}\bigr\}_{\nu} \in \mathscr{K}$ for the deterministic dynamical system of Equation}~\eqref{Eq5}); {while minimizing the asymptotic exit rate with which the controlled-diffusion process $x^{\epsilon}(t)$ exits from the given domain $D$.}

{Note that a closer look at Proposition~}\ref{P1} {suggests an optimization problem that links between the minimum asymptotic exit rate $r(\mathcal{K}^{\ast})$ for the controlled-diffusion process $x^{\epsilon}(t)$ and a certain $m$-tuple of (sub)-optimal admissible linear feedback operators $\mathcal{K}^{\ast} = (\mathcal{K}_1^{\ast}, \mathcal{K}_2^{\ast}, \ldots \mathcal{K}_m^{\ast}) \in \mathscr{K}$ which corresponds to the deterministic dynamical system of Equation~}\eqref{Eq5}. {Namely, if the statement in the above proposition holds true, then there exists at least one $m$-tuple of linear feedback operators $\mathcal{K}^{\ast} \in \mathscr{K}$ such that}
\begin{align}
\mathcal{K}^{\ast} \in \argmin_{\mathcal{K} \in \mathscr{K}} \Biggl \{\limsup_{T \rightarrow \infty} \,\inf_{\substack{\varphi(t) \in C_{0T}([0,T], \mathbb{R}^d)\\ \varphi(0)=x_0}} \frac{1}{T} \biggl \{ S_{0T}^{\mathcal{K}}(\varphi(t)) \, \Bigl \vert \,\varphi(t) \in D \cup \partial D, \forall t \in [0, T] \biggr\} \Biggr\}. \label{Eq37}
\end{align}
{Moreover, if such a solution exists, then the admissible control $u_i^{\ast}(\cdot)=\bigl(\mathcal{K}_i^{\ast} x^{\epsilon}\bigr)(\cdot) \in \mathcal{U}_i$, $\forall t \ge 0$, with $\mathcal{K}_i^{\ast} \in \mathscr{K}_i[\mathcal{X}, \mathcal{U}_i]$, for $i = 1, 2, \ldots, m$, is essentially a maximizing measurable selector of the following Hamilton-Jacobi-Bellman equation related to an optimal control problem}
\begin{align}
\left.\begin{array}{c}
 \max_{u}  \bigl \{ \mathcal{L}_{\epsilon}\, \upsilon_{(\epsilon, x_0)}^{\ast}(x, u) + \lambda^{\ast} \,\upsilon_{(\epsilon, x_0)}^{\ast}(x) \bigr\}, \quad \forall x \in D  \\
\quad\quad~~ ~ \upsilon_{(\epsilon, x_0)}^{\ast}(x) = 0, \quad \forall x \in \partial D \
\end{array}\right\}, \label{Eq38}
\end{align}
{where }
\begin{align*}
 \mathcal{L}_{\epsilon}(\cdot)(x, u) = \Bigl \langle \bigtriangledown(\cdot), \Bigl(A x + \sum\nolimits_{i=1}^m B_i u_i \Bigr) \Bigr\rangle + \frac{\epsilon}{2} \operatorname{tr}\Bigl \{\sigma(x)\sigma^T(x)\bigtriangledown^2(\cdot) \Bigr\},
\end{align*}
{$\upsilon_{(\epsilon, x_0)}^{\ast} \in C^2(D) \cap C(D \cup \partial D)$, with $\upsilon_{(\epsilon, x_0)}^{\ast} > 0$ in $D$, and $u(\cdot) = (u_1(\cdot), u_2(\cdot), \ldots, u_m(\cdot)) \in \prod\nolimits_{i=1}^m \mathcal{U}_i$ (e.g., see} \cite[Theorem~1.4(a)]{QuaSi08} {for additional discussions; and cf.} \cite{Fle78} or \cite{EvaIsh85}).

\begin{remark} \label{R3}
{Note that if $u^{\ast}(\cdot)$ is the maximizing measurable selector for $\argmax \bigl \{ \mathcal{L}_{\epsilon} \upsilon_{(\epsilon, x_0)}^{\ast}(x, \cdot) \bigr\}$, $x \in D$, then the principal eigenvalue is given by}
\begin{align*}
\lambda_{\epsilon}^{\mathcal{K}^{\ast}} = -\limsup_{T \rightarrow \infty} \frac{1}{T} \log \mathbb{P}_{\epsilon}^{\mathcal{K}^{\ast}} \bigl\{ \tau_{D}^{\epsilon} > T\bigr\}, 
\end{align*}
{where the probability distribution $\mathbb{P}_{\epsilon}^{\mathcal{K}^{\ast}}\bigl\{\cdot\bigr\}$ is conditioned with respect to the admissible control $u^{\ast}$ and the initial point $x_0 \in D$. Furthermore, if $\lim_{\epsilon \rightarrow 0} \lambda_{\epsilon}^{\mathcal{K}^{\ast}} < \infty$, for some $x_0 \in D$, then the maximum closed invariant set $\Lambda_D^{\mathcal{K}^{\ast}} \subset D \cup \partial D$ for the deterministic dynamical system (under the action of the $m$-tuple of linear feedback operators $\mathcal{K}^{\ast} \in \mathscr{K}$) is nonempty (cf. Remark}~\ref{R1}). 
\end{remark}

\begin{remark} \label{R4}
{Finally, it is worth remarking that Proposition~}\ref{P1} {is useful for selecting the most appropriate $m$-tuple of admissible linear feedback operators from the set $\bigl\{\mathcal{K}^{\ast}\bigr\}_{\nu} \in \mathscr{K}$ that confines the controlled-diffusion process $x^{\epsilon}(t)$ to the prescribed domain $D$ for a longer duration.} 
\end{remark}

\end{document}